\documentclass[11pt,twoside,a4paper]{article}
\usepackage{mathrsfs}
\usepackage{mathrsfs}
\usepackage{mathrsfs}
\usepackage{mathrsfs}
\usepackage{mathrsfs}
\usepackage{mathrsfs}
\usepackage{mathrsfs}
\usepackage{stmaryrd}
\usepackage{amsfonts}
\usepackage{amsmath,amssymb}
\usepackage{mathrsfs}
\usepackage{hyperref}
\usepackage{graphicx}
\usepackage{graphicx}

\newtheorem{thm}{Theorem}[section]
\newtheorem{cor}[thm]{Corollary}

\newtheorem{prop}[thm]{Proposition}

\newtheorem{rem}[thm]{Remark}

\newtheorem{que}[thm]{Question}
\numberwithin{equation}{section}\allowdisplaybreaks

\begin{document}
 
\title{\bf Modulation  Spaces  and  Nonlinear Evolution Equations}

\author{Michael Ruzhansky$^a$, \ Mitsuru Sugimoto$^b$, \ Baoxiang Wang$^c$\footnote{
$^a$Department of Mathemaics,
 Imperial College London,
180 Queens Gate, London SW7 2AZ
 United Kingdom, \
Email: m.ruzhansky@imperial.ac.uk, \  $^b$Graduate School of Mathematics,
Nagoya University,
Furocho, Chikusa-ku,
Nagoya 464-8602 Japan, Email: sugimoto@math.nagoya-u.ac.jp, \ $^c$LMAM, School of Mathematical Sciences,
Peking University, and Beijing International Center of Mathematical Research,
Beijing 100871,
PR of China, Email:wbx@math.pku.edu.cn}}

\maketitle

\begin{abstract}
We survey some recent progress on modulation spaces and the well-posedness results for a class of nonlinear evolution equations by using the frequency-uniform localization techniques.
\end{abstract}

\section{Definitions of modulation spaces}

The modulation space was introduced by Feichtinger \cite{Fei83} in 1983. We denote by $\mathfrak{S}$ the Schwartz space and by $\mathfrak{S}'$ its dual space. Let $f,g $ be Schwartz functions, and set
\begin{equation*}
V_{g}f(x,\omega)=\int_{\mathbb{R}^{n}}e^{-{\rm i} t
\omega}\overline{g(t-x)}f(t)dt.
\end{equation*}
$V_g f$ is said to be the short-time Fourier transform of $f$. In the following, if there is no explanation, we will always assume that
$$
-\infty< s<\infty, \ \  0<p,q \le \infty, \ \ \langle \cdot\rangle =(1+|\cdot|^2)^{1/2}.
$$
The norm on modulation spaces is given by
\begin{equation} \label{moddef}
\|f\|^{\rm o}_{M^{s}_{p,q}}=\left(\int_{\mathbb{R}^{n}}\left(\int_{\mathbb{R}^{n}}|V_{g}f(x,\omega)|^{p}dx\right)^{q/p}\langle
\omega\rangle^{sq}d\omega\right)^{1/q},
\end{equation}
with a natural modification for $p,q=\infty$.
Applying the frequency-uniform localization techniques, one can give an equivalent norm on modulation spaces. Let $Q_k$ be the unit cube with its center at $k$,   $\{Q_k\}_{k\in \mathbb{Z}^n}$ constitutes a decomposition of  $\mathbb{R}^n$. Such kind of decomposition  goes back to the work of N. Wiener \cite{Wie32}, and we  say that it is the Wiener decomposition of  $\mathbb{R}^n$.
We can roughly write
\begin{equation} \label{fud}
\Box_k \sim \mathfrak{F}^{-1} \chi_{Q_k} \mathfrak{F}, \quad k\in \mathbb{Z}^n,
\end{equation}
where  $\chi_E$ denotes the characteristic function on the set $E$, and
$\mathfrak{F}$ is the Fourier transform. Since $Q_k$ is a translation of $Q_0$, $\Box_k$ ($k\in \mathbb{Z}^n$) have the same localized structures in the frequency space, which are said to be the {\it frequency-uniform decomposition operators}. Similar to Besov spaces, one can use $\{\Box_k\}_{k\in \mathbb{Z}^n}$ and $\ell^q(L^p)$ to generate a class of function spaces, so called modulation spaces for which the norm is defined by
$$
\|f\|_{M^s_{p,q}} = \left(\sum_{k \in \mathbb{Z}^n} \langle
k\rangle^{sq} \|\,\Box_k f\|_p^q \right)^{1/q},
$$
with the natural modifications for $p,q=\infty$.
We now give an exact definition on frequency-uniform decomposition operators. Since $\chi_{Q_k}$ can not be differentiated, it is convenient to replace $\chi_{Q_k}$ in \eqref{fud} by a smooth cut-off function.  Let $\rho\in
\mathfrak{S}(\mathbb{R}^n)$, $\rho:\, \mathbb{R}^n\to [0,1]$ be a smooth function verifying $\rho(\xi)=1$ for $|\xi|_\infty \le 1/2$ and $\rho(\xi)=0$ for $|\xi|_\infty \ge 1$\footnote{For $\xi=(\xi_1,...,\xi_n)$, $|\xi|_\infty:=\max_{i=1,...,n}|\xi_i|$.}. Let $\rho_k$ be a translation of $\rho$,
\begin{equation}
\rho_k (\xi) = \rho (\xi- k), \quad k\in \mathbb{Z}^n. \label{freq2.1}
\end{equation}
We see that $\rho_k (\xi)=1$ in $Q_k$ and so, $\sum_{k\in \mathbb{Z}^n}\rho_k(\xi)\ge 1$
for all $\xi\in \mathbb{R}^n$. Denote
\begin{equation}
\sigma_k (\xi)= \rho_k(\xi) \left(\sum_{\ell\in \mathbb{
Z}^n}\rho_\ell(\xi)\right)^{-1}, \quad k\in \mathbb{Z}^n. \label{freq2.2}
\end{equation}
Then we have
\begin{equation}
 \left\{\begin{array}{l}
|\sigma_k(\xi)| \ge c, \quad \forall \; \xi \in Q_k,\\
 {\rm supp}\, \sigma_k \subset \{\xi: |\xi-k|_\infty \le 1\},\\
\sum_{k\in \mathbb{Z}^n} \sigma_k (\xi) \equiv 1, \quad \forall \;
\xi \in
\mathbb{R}^n,\\
|D^\alpha\sigma_k(\xi)| \le C_{|\alpha|}, \quad \forall \; \xi \in \mathbb{
R}^n, \; \alpha \in (\mathbb{N} \cup \{0\})^n
.
\end{array}\right. \label{freq2.3}
\end{equation}
Hence, the set
\begin{equation}
\Upsilon_n = \{ \{\sigma_k\}_{k\in \mathbb{Z}^n}: \;
\{\sigma_k\}_{k\in \mathbb{Z}^n}\;\; \mbox {satisfies  \eqref{freq2.3}}
\} \label{freq2.4}
\end{equation}
is non-void.  If there is no confusion, in the sequel we will write $\Upsilon=\Upsilon_n$.
Let $\{\sigma_k\}_{k\in \mathbb{Z}^n}\in \Upsilon$.
The frequency-uniform decomposition operators can be exactly defined by
\begin{equation}
\Box_k := \mathfrak{F}^{-1} \sigma_k \mathfrak{F}, \quad k\in \mathbb{
Z}^n. \label{freq2.5}
\end{equation}
For simplicity, we write
$M^0_{p,q}=M_{p,q}$. The space $M^s_{p,q}$  is said to be the modulation space.
One can prove that $\|\cdot\|^{\rm o}_{M^{s}_{p,q}}$ and  $\|\cdot\|_{M^{s}_{p,q}}$ are equivalent norms; cf. \cite{Fei83} for a proof on modulation spaces defined on Abelian groups, and \cite{WaHe07} for a straightforward proof.

\section{Questions on the Schr\"odinger equation  \label{WhyFUD}}

Let $S(t)=e^{{\rm i}t\Delta}$ be the Schr\"odinger semi-group. It is known that  $S(t)=e^{{\rm i}t\Delta}: L^p \to L^p$
is bounded if and only if $p=2$. Indeed, one can show that
\begin{equation}
 \sup_{\varphi \in  \mathfrak{S}\setminus\{0\}} \frac{\|S(t) \varphi \|_{L^p ({\Bbb R}^n)}}{\|\varphi\|_{L^{p}({\Bbb R}^n)}}
=\infty, \quad p\not= 2,   \label{Ppinfty}
\end{equation}
by taking $\varphi= e^{-(a+{\rm i}b)|x|^2}$. So, a natural question is the following

\begin{que}\label{Q21}
Is there any Banach function space $X$ satisfying
\begin{equation}
  \sup_{\varphi\in X \setminus \{0\}} \frac{\|S(t) \varphi \|_{X({\Bbb R}^n)}}{\|\varphi\|_{X({\Bbb R}^n)}}
< \infty \, ?\label{Ppfty}
\end{equation}
\end{que}
Recall that (cf. \cite{CW, Cazebook03})
\begin{equation} \label{unmod}
\| S(t) f\|_p \le C  |t|^{-n(1/2-1/p)}\| f\|_{p'},
\end{equation}
where $2\le p\le \infty, \ 1/p+1/p'=1$. Estimate \eqref{unmod} contains singularity at $t=0$, and another natural question is
\begin{que}\label{Q22}
Is there any Banach function space $X$ satisfying the following truncated decay
\begin{equation}
   \|S(t) \varphi \|_{X({\Bbb R}^n)}  \le C  (1+|t|)^{-\delta(X)} \|\varphi\|_{X'({\Bbb R}^n)} ?
 \label{Ppfty2}
\end{equation}
where $\delta>0$, $X'$ denotes the dual space of $X$, $X$ and $X'$ have the same regularity.
\end{que}

Consider the Cauchy problem for the nonlinear Schr\"odinger equation (NLS)
\begin{equation}
{\rm i}u_t + \Delta u =   |u|^{2\kappa} u, \quad u(0,x)= u_0(x), \  \kappa \in \mathbb{N}. \label{NLSmoda}
\end{equation}
It is known that if $u$ is a solution of \eqref{NLSmoda}, so is $u_\mu (t,x)= \mu^{1/\kappa} u(\mu^2 t, \mu x)$ with the initial datum $\mu^{1/\kappa} u_0(\mu x)$.
Taking notice of
$$
\|u_\mu (0)\|_{\dot H^s} = \mu^{s-n/2+ 1/\kappa} \|u_0\|_{\dot H^s}, \ \ \mu>0,
$$
one sees that $s=s_\kappa : = n/2- 1/\kappa$ is the index such that $\|u_\mu (0)\|_{\dot H^s}$ is invariant for all $\mu>0$. $\dot H^{s_\kappa}$ is said to be a critical space for NLS \eqref{NLSmoda}; cf. \cite{Cazebook03}. Up to now, we can solve \eqref{NLSmoda} in $H^s$  for $s\ge s_\kappa \ge 0$, which correspond to the critical and sub-critical cases in $H^s$. However, we cannot solve  \eqref{NLSmoda} in $H^s$ in the case $s<s_\kappa$, where $H^s$ is said to be the super-critical space. Our question is
\begin{que}\label{Q23}
Let $0\le s<s_\kappa$. Are there any initial data $u_0\in H^s$ so that NLS \eqref{NLSmoda} is still locally and globally well posed?
\end{que}

In order to answer the above questions, modulation spaces play important roles.  Roughly speaking, Schr\"odinger semi-group is bounded in modulation spaces and it satisfies a truncated decay in modulation spaces. Moreover, the modulation space  $M_{2,1}$ enjoys lower derivative regularity and we can solve NLS \eqref{NLSmoda} in $M_{2,1}$ for all $\kappa\in \mathbb{N}$. $M_{2,1}$ contains a class of data which are out of the control of $H^{s_\kappa}$ if $s_\kappa >0$.

Finally, we mention that all the questions above make sense for other equations as well.
For example, the propagators $W_\pm(t)=e^{\pm{\rm i} t\sqrt{-\Delta}}$ for the wave equation
display the loss of regularity for $p\not=2$. More precisely, we have that for $1<p<\infty$ and $t>0$, the operator
$W_+(t):W^{s,p}(\mathbb R^n)\to L^p(\mathbb R^n)$ is (locally) bounded if and only if
$s\geq (n-1)|1/p-1/2|$. The same holds for $W_-(t)$, and Question \ref{Q21} is valid.
Also, a variant \eqref{unmod} holds\footnote{with $n$ replaced by $n-1$.}, so
 Question \ref{Q22} is also valid, as well as the corresponding version of
 Question \ref{Q23}.

We note that the dependence of estimates and large time asymptotics in
Sobolev spaces for the propagator
$e^{{\rm i}t a(D)}$ on $a(\xi)$ is known, especially in the first order case.
For example, let $a\in C^\infty({\mathbb R}^n\backslash 0)$ be a real-valued positively homogeneous
function, satisfying $\nabla a(\xi)\not=0$ for all $\xi\not=0$. If we denote
$$
 k:=\max_{|\xi|=1} {\rm rank} \nabla_\xi^2 a(\xi),
$$
then there is a loss of derivatives in $L^p$ depending on $k$. More
precisely, as a special case of estimates for
Fourier integral operators for $1<p<\infty$, it was shown in \cite{Ru99} (see also \cite{Ru00}) that
if $e^{{\rm i}a(D)}:W^{s,p}\to L^p$ is locally bounded then
$s\geq k|1/p-1/2|$. The dispersive estimates\footnote{or rather the value of $\kappa$.}
\begin{equation*} \label{unmoda}
\| e^{{\rm i}t a(D)} f\|_p \le C  |t|^{-\kappa(1/2-1/p)}\| f\|_{p'},
\end{equation*}
where $2\le p\le \infty, \ 1/p+1/p'=1$, depend on the
convexity of the level sets $\Sigma=\{\xi: a(\xi)=1\}$ as well as
on the maximal orders of contacts between $\Sigma$ and its
tangent lines. For a comprehensive analysis of such types of
estimates and their applications to the Strichartz estimates we
refer to \cite{RuSm10}.

\section{Some results on modulation spaces}

Roughly speaking,  frequency-uniform decomposition operators combined with function spaces $\ell^q(L^p)$ produce modulation spaces.
During the past thirty years, the importance of the frequency-uniform decomposition in applications seems to be not
fully recognized,  and it is even not mentioned in Gr\"ochenig's famous book \cite{Groch01}. However,
from PDE point of view, the combination of frequency-uniform decomposition operators and  Banach function spaces $\ell^q(X(\mathbb{R}^n))$\footnote{$X$ is a Banach function space defined in $\mathbb{R}^n$.}
 seems to be important in making nonlinear estimates, which contains an automatic decomposition on high-low frequencies.

\subsection{Basic properties on modulation spaces}

 \begin{prop}[Completeness]\label{ModP1}
Let $0<p,q\le \infty$ and $s\in \mathbb{R}$.
\begin{itemize}
\item[\rm (1)] $M^s_{p,q}$ is a (quasi-) Banach space. Moreover, if
 $1\le p,q\le \infty$, then $M^s_{p,q}$  is a Banach space.
\item[\rm (2)]  $\mathfrak{S}(\mathbb{R}^n)\subset M^s_{p,q} \subset
\mathfrak{S}'(\mathbb{R}^n)$.
\item[\rm (3)] Let $0<p,q<\infty$, then $\mathfrak{S}(\mathbb{R}^n)$  is dense in $M^s_{p,q}$.
\end{itemize}
\end{prop}

\begin{prop}[Equivalent norm]\label{ModP2}
Let
$\{\sigma_k\}_{k\in {\Bbb Z}^n},\  \{\varphi_k\}_{k\in {\Bbb
Z}^n}\in \Upsilon$. Then $\{\sigma_k\}_{k\in {\Bbb Z}^n}$ and  $
\{\varphi_k\}_{k\in {\Bbb Z}^n}$ generate equivalent norms on $M^s_{p,q}$.
\end{prop}

\medskip

Proposition \ref{ModP2} indicates that one can choose
$\{\sigma_k\}_{k\in \mathbb{Z}^n}\in \Upsilon_n$ according to our requirement. In applications of PDE,  it is convenient to use the following $\{\sigma_k\}_{k\in \mathbb{Z}^n}\in \Upsilon_n$.
Let $\{\eta_k\}_{k\in \mathbb{Z}}\in \Upsilon_1$, we denote
\begin{equation}
\sigma_k(\xi): = \eta_{k_1}(\xi_1)...\eta_{k_n}(\xi_n), \label{freq-u-dec}
\end{equation}
then we have $\{\sigma_k\}_{k\in \mathbb{Z}^n}\in \Upsilon_n$. the above $\sigma_k(\xi)$ realizes the separation of different variables.

\begin{prop}[Embedding]\label{ModP3}
Let $s_1, s_2 \in \mathbb{R}$ and $0<p_1, p_2, q_1, q_2\le \infty$.
\begin{itemize}
\item[\rm (1)] If $s_2 \le s_1$,  $p_1\le p_2$ and $q_1\le q_2,$ then
$M^{s_1}_{p_1,q_1}  \subset M^{s_2}_{p_2,q_2}. $
\item[\rm (2)] If $q_2< q_1$ and
 $s_1-s_2>n/q_2-n/q_1$, then $M^{s_1}_{p,q_1}\subset M^{s_2}_{p,q_2}$.
\end{itemize}
\end{prop}

\begin{prop}[Dual space]\label{ModP4}
Let $s\in \mathbb{R}$ and $0<p,q<\infty$. If $p\ge 1$, we denote $1/p+1/p'=1$; If $0<p<1$, we write $p'=\infty.$ Then
\begin{equation}
(M^s_{p,q})^*= M^{-s}_{p',q'}.
 \label{mod-dual}
\end{equation}
\end{prop}

If $p\ge 1$, Proposition \ref{ModP4} is similar to that of Besov spaces, however, if $0<p<1$, the result is quite different from that of Besov spaces.    The details of the proof of Proposition \ref{ModP4} can be found in \cite{WaHe07} by following the proof of the relevant result in Besov spaces.

\begin{rem}
If $p,q\in [1,\infty]$, Propositions \ref{ModP1} and \ref{ModP4}
were obtained by Feichtinger \cite{Fei83}. In \cite{WaHe07,Wang8}, the cases $0<p<1$ and $0<q<1$ were considered.

Soon after the work \cite{Wang8},  Kobayashi \cite{Kub06} independently defined $M_{p,q}$ for all $0<p,q\le \infty$ and obtained Proposition \ref{ModP1}.
Almost at the same time as \cite{WaHe07}, Kobayashi \cite{Kub062} discussed the dual space of $M_{p,q}$ and obtained partial results of Proposition \ref{ModP4}:
if $0<p<1$ or $1<q< \infty$, he obtained $M_{p',q'} \subset (M_{p,q})^* \subset M_{\infty, \infty}$.
For the other cases, he showed $(M_{p,q})^*= M_{p',q'}.$  Recently, by using the molecular decomposition techniques of modulation spaces, Kobayashi and Sawano \cite{KoSa10} reconsidered the dual space of $M^s_{p,q}$ and they also obtained the result of  Proposition \ref{ModP4}. It is worth to mention that Triebel \cite{Tr83} introduced a class of generalized modulation spaces for all indices $0<p,q\le \infty$, however, those spaces have no complete norms, which seems harder to use in the study of PDEs.
\end{rem}

\subsection{Inclusions between Besov and modulation spaces}

From the definitions, we see that  Besov spaces and modulation spaces are rather similar, both of them are the combinations of frequency decomposition operators and function spaces $\ell^q(L^p)$. In fact, we have the following inclusion results.

\begin{thm}[Embedding] \label{embedding1}
Let $0< p, q \le \infty$ and $s_1, s_2 \in \mathbb{R}$. We have the following results.

{\rm (1)} $ B^{s_1 }_{p,q} \subset M^{s_2}_{p,q}$  if and only if
$s_1\ge s_2 + \tau (p,q)$, where
\begin{equation}
\tau (p,q)= \max \left\{0, \; n\left(\frac{1}{q}-\frac{1}{p}\right),
\; n\left(\frac{1}{q}+\frac{1}{p}-1\right) \right\}; \nonumber
\end{equation}

{\rm (2)} $ M^{s_1 }_{p,q} \subset B^{s_2}_{p,q}$ if and only if
$s_1\ge s_2 + \sigma (p,q)$, where
\begin{equation}
\sigma (p,q)= \max \left\{0, \;
n\left(\frac{1}{p}-\frac{1}{q}\right), \;
n\left(1-\frac{1}{p}-\frac{1}{q} \right) \right\}. \nonumber
\end{equation}
\end{thm}

%\begin{figure}
% \centering
% \includegraphics[width=8cm,height=8cm]{fig1}
% \begin{center}
% \begin{minipage}{11cm}
% \caption{\footnotesize The distribution of $\tau(p,q)$
% $\mathbb{R}^2_{+}$:  $\tau(p,q)=n(\frac{1}{q}-\frac{1}{p})$ in $S_1$;
% $\tau(p,q)= 0$ in $S_2$; $\tau(p,q)=n(\frac{1}{p}+\frac{1}{q}-1)$ in $S_3$.}
% \end{minipage}
% \end{center}
% \end{figure}

% \begin{figure}
% \centering
% \includegraphics[width=8cm,height=8cm]{fig2}
% \begin{center}
% \begin{minipage}{11cm}
% \caption{\footnotesize The distribution of $\sigma(p,q)$: $\sigma(p,q)= 0$ in $R_1$; $\sigma(p,q)=n(\frac{1}{p}-\frac{1}{q})$ in $R_2$; $\sigma(p,q)=n(1-\frac{1}{p}-\frac{1}{q})$ in $R_3$. }
% \end{minipage}
% \end{center}
% \end{figure}

The inclusions  between Besov and modulation spaces in the cases $(1/p, 1/q)\in [0,1]^2$ were first discussed by Gr\"obner \cite{Grob92} and he has never published his results. When $(1/p, 1/q)$ is in the vertices of the square $[0,1]^2$, Gr\"obner's results are optimal.  Afterwards, Toft \cite{Toft04} obtained the sufficiency of Theorem \ref{embedding1} in the cases $(1/p, 1/q)\in [0,1]^2$.  Sugimoto and  Tomita \cite{SuTo07} showed  the necessity of the first inclusion of Theorem \ref{embedding1} in the cases $(1/p, 1/q)\in [0,1]^2$, and  by duality the second inclusion is also sharp if  $(1/p, 1/q)\in [0,1]^2$ and $p,q\not=\infty$.  Sugimoto and Tomita's idea is to use Feichtinger's norm and the dilation property of modulation spaces. In \cite{WaHu07,WaHe07,Wang8} the authors proved the conclusions of Theorem \ref{embedding1} by using frequency-uniform decomposition techniques.

\begin{cor} \label{embedding2}
 We have the following inclusions.
$$
B^{s+n/2}_{2,1} \subset M^s_{2,1} \subset B^s_{2,1}, \quad B^{s+n}_{\infty,1} \subset M^s_{\infty,1} \subset B^s_{\infty,1}.
$$
\end{cor}

The above embedding theorem is of importance for the study of nonlinear PDEs.
As for the inclusions between $L^p$-Sobolev spaces and modulation spaces,
it has been explicitly determined by a recent work of Sugimoto and Kobayashi
\cite{KoSu11}:

\begin{thm}\label{embedding3}
Let $1 \leq p,q \leq \infty$ and and $s_1, s_2 \in \mathbb{R}$.
Then $ H^{s_1}_p \subset M_{p,q}^{s_2}$
if and only if one of the following conditions is satisfied$:$

{\rm (1)} $q \geq  p> 1$, $s_1 \geq s_2+\tau(p,q);$
\quad\, \hspace{5.5mm}
{\rm (2)} $p>q$, $s_1> s_2+\tau(p,q);$

{\rm (3)} $p=1, q=\infty$, $s_1 \geq s_2+\tau(1,\infty);$
\quad\,
{\rm (4)} $p=1, q \not=\infty$, $s_1> s_2+\tau(1,q)$,

\noindent
and $M_{p,q}^{s_1} \subset H^{s_2}_p$
if and only if one of the following conditions is satisfied$:$

{\rm (1)} $q \leq p< \infty$, $s_1 \geq s_2+ \sigma(p,q);$
\qquad\,\,\,
{\rm (2)} $p<q$, $s_1 > s_2+ \sigma(p,q);$

{\rm (3)} $p=\infty, q=1$, $s_1 \geq s_2+\sigma(\infty,1);$
\quad\,\,\,
{\rm (4)} $p=\infty, q \not=1$, $s_1> s_2+ \sigma(\infty,q)$,

\noindent
where $\tau(p,q)$ and $\sigma(p,q)$ are the same indices as in Theorem
\ref{embedding1}.
\end{thm}

\subsection{Dilation property of modulation spaces}

Roughly speaking, homogeneous Sobolev spaces and their generalizations including homogeneous Besov and Triebel spaces have scales like
$$
\|f(\lambda \ \cdot)\|_X \sim \lambda^{\theta (X)} \|f\|_X.
$$
However, the scaling properties of modulation spaces  are  very complicated, which is quite different from the classical Sobolev spaces. The following result is due to Sugimoto and Tomita. Let  $\tau(p,q)$ and $\sigma(p,q)$ be as in Theorem \ref{embedding1}.  Denote
$$
\mu_1(p,q) = \tau(p,q) -n/p, \ \  \mu_2(p,q) = -\sigma(p,q) -n/p
$$
\begin{thm}  [\cite{SuTo07}]
Let $1\le p,q\le \infty$. The following are true.
\begin{itemize}
\item[(1)] There exists a constant $C>0$ such that
$$
C^{-1} \lambda^{\mu_2(p,q)} \|f\|_{M_{p,q}} \le \|f(\lambda \ \cdot)\|_{M_{p,q}} \le  C  \lambda^{\mu_1(p,q)} \|f\|_{M_{p,q}}
$$
for all $f\in M_{p,q}$ and $\lambda\ge 1$. Conversely, if there is a constant $C>0$ such that
$$
C^{-1} \lambda^{\beta} \|f\|_{M_{p,q}} \le \|f(\lambda \ \cdot)\|_{M_{p,q}} \le  C  \lambda^{\alpha} \|f\|_{M_{p,q}}
$$
for all $f\in M_{p,q}$ and $\lambda\ge 1$, then $\alpha \ge \mu_1(p,q)$ and $\beta \le \mu_2(p,q)$.

\item[(2)] There exists a constant $C>0$ such that
$$
C^{-1} \lambda^{\mu_1(p,q)} \|f\|_{M_{p,q}} \le \|f(\lambda \ \cdot)\|_{M_{p,q}} \le  C  \lambda^{\mu_2(p,q)} \|f\|_{M_{p,q}}
$$
for all $f\in M_{p,q}$ and $0<\lambda\le 1$. Conversely, if there is a constant $C>0$ such that
$$
C^{-1} \lambda^{\alpha} \|f\|_{M_{p,q}} \le \|f(\lambda \ \cdot)\|_{M_{p,q}} \le  C  \lambda^{\beta} \|f\|_{M_{p,q}}
$$
for all $f\in M_{p,q}$ and $0< \lambda\le 1$, then $\alpha \ge \mu_1(p,q)$ and $\beta \le \mu_2(p,q)$.
\end{itemize}

\end{thm}

\section{NLS and NLKG in modulation spaces}

As indicated in \S \ref{WhyFUD},  the dispersive semi-group combined with the frequency-uniform decomposition operator has some advantages and we discuss them in this section. The results of this section can be found in \cite{BenOk08,BGOR,CorNik09,Wang8,WaHe07}.

\subsection{Schr\"odinger and Klein-Gordon semigroup in modulation spaces}

Let $S(t)= e^{{\rm i}t \triangle}$ denote the Schr\"odinger semi-group. In \cite{Wang8}, Wang, Zhao and Guo obtained the uniform boundedness for the Ginzburg-Landau semi-group $L(t)= e^{(a+{\rm i})t \triangle}$ ($a>0$) in modulation spaces and their proof is also adapted to the Schr\"odinger semi-group ($a=0$ in $L(t)$).

\begin{prop}[Uniform boundedness of $S(t)$ in $M^s_{p,q}$ ]\label{propmodpp.1}
Let $s\in \mathbb{R}$, $1\le p \le \infty$ and  $0<q<\infty$.
Then we have
\begin{equation}
\|S(t)f\|_{M^s_{p,q}} \le C  (1+|t|)^{n|1/2-1/p|} \|
f\|_{M^s_{p,q}}. \label{modpp.5}
\end{equation}
\end{prop}

Shortly after the work \cite{Wang8}, Proposition \ref{propmodpp.1} is independently obtained by B\'{e}nyi,
Gr\"ochenig, Okoudjou and Rogers in \cite{BGOR} and their result  contains more general semi-group $e^{{\rm i}t(-\Delta)^\alpha}$ with $\alpha\le 1$, whose proof is based on the short-time frequency analysis technique.   Miyachi,   Nicola,   Riveti,  Taracco and Tomita \cite{MNRTT09} were able to consider the case $\alpha>1$, Chen, Fan and Sun \cite{CFS10} obtained some refined estimates for $e^{{\rm i}t(-\Delta)^\alpha}$ with any $\alpha>0$.

Now we consider the truncated decay of $S(t)$.
\begin{prop}\label{propmodecay.1}
Let $s\in \mathbb{R}$, $2\le p<\infty$, $1/p+1/p'=1$ and $0<q<\infty$.
Then we have
\begin{equation}
\|S(t)f\|_{M^s_{p,q}} \le C  (1+|t|)^{-n(1/2-1/p)} \|
f\|_{M^s_{p',q}}. \label{modecay.4}
\end{equation}
\end{prop}

Propositions \ref{propmodpp.1} and \ref{propmodecay.1} answer Questions 2.1 and 2.2, and moreover, they are optimal in the sense that the powers of time variable are sharp, cf. \cite{CorNik092}. Now we consider the truncated decay estimate for the Klein-Gordon semi-group $G(t)= e^{{\rm i}t \omega^{1/2}}$ where $\omega=I-\Delta$.

\begin{prop}\label{propmodpp.g1}
Let $s\in \mathbb{R}$, $1\le p \le \infty$ and  $0<q<\infty$.
Then we have
\begin{equation}
\|G(t)f\|_{M^s_{p,q}} \le C  (1+|t|)^{n|1/2-1/p|} \|
f\|_{M^s_{p,q}}. \label{modpp.5_2}
\end{equation}
\end{prop}

It is known that
 $G(t)$ satisfies the following
  $L^p-L^{p'}$ estimate
\begin{equation}
\|G(t)f\|_{H^{-2\sigma(p)}_p} \le C  |t|^{-n(1/2-1/p)} \|
f\|_{p'}, \label{modecay.5}
\end{equation}
where
\begin{equation}
2\le p<\infty, \quad 2\sigma(p) = (n+2)\left(\frac{1}{2}-
\frac{1}{p}\right). \label{modecay.6}
\end{equation}
From \eqref{modecay.5} it follows that

\begin{prop}\label{propmodecay.2}
Let $s\in \mathbb{R}$, $2\le p<\infty, \, 1/p+1/p'=1$, $0<q<\infty$,
$\theta\in [0,1]$ and $\sigma(p)$ is as in
 \eqref{modecay.6}.  Then we have
\begin{equation}
\|G(t)f\|_{M^s_{p,q}} \le C  (1+|t|)^{-n\theta (1/2-1/p)} \|
f\|_{M^{s+2\sigma(p)\theta}_{p',q}}. \label{modecay.14}
\end{equation}
\end{prop}

\subsection{Strichartz estimates in modulation spaces}

For convenience, we write
\begin{equation}
\| f\|_{\ell^{s, q}_\Box (L^\gamma(I, L^p)) }
= \left(\sum_{k\in \mathbb{Z}^n} \langle k\rangle^{sq}\|\Box_k f\|^q_{L^\gamma(I, L^p)}\right)^{1/q}, \label{fsud}
\end{equation}
$\ell^{q}_\Box (L^\gamma(I, L^p)):= \ell^{0, q}_\Box (L^\gamma(I, L^p))$, $\ell^{q}_\Box (L^p_{x,t\in I}):= \ell^{q}_\Box (L^p(I, L^p))$.
Recall that the truncated decay can be generalized to the following estimate
\begin{equation}
\|U(t)f\|_{M^\alpha_{ p,q}} \le C  (1+|t|)^{-\delta} \| f\|_{M_{
p',q}}, \label{modstr.1}
\end{equation}
where $2\le p<\infty$, $1\le q< \infty$, $\alpha=\alpha(p)\in
\mathbb{R}$, $\delta=\delta(p)>0$, $\alpha$ and $ \delta$ are independent of
 $t\in \mathbb{R}$,   $U(t)$ is a dispersive semi-group,
\begin{equation}
U(t)= \mathfrak{F}^{-1} e^{{\rm i}t P(\xi)} \mathfrak{F}, \label{modstr.2}
\end{equation}
and $P(\cdot): \mathbb{R}^n\to \mathbb{R}$ is a real valued function.
In the sequel we will assume that $U(t)$ satisfies conditions \eqref{modstr.1} and \eqref{modstr.2}, from which we can get some Strichartz inequalities for $U(t)$ in modulation spaces (cf. \cite{WaHe07}).

\begin{prop}[Strichartz inequalities]\label{propmodstr5.1}
Let $U(t)$ satisfy \eqref{modstr.1} and \eqref{modstr.2}.  For any
$\gamma\ge 2 \vee (2/\delta)$, we have
\begin{equation}
\|U(t)f\|_{\ell^{\alpha/2, q}_\Box (L^\gamma(\mathbb{R}, L^p)) }
\le C   \| f\|_{M_{2,q}}. \label{modstr.3}
\end{equation}
In addition, if $\gamma\ge q$, then we have
\begin{equation}
\|U(t)f\|_{L^\gamma(\mathbb{R}, M^{\alpha/2}_{p,q})} \le C  \|
f\|_{M_{2,q}}. \label{modstr.4}
\end{equation}
\end{prop}

Denote
\begin{equation}
(\mathfrak{U} f)(t)= \int^t_0  U(t-s)f(s,\cdot) ds. \label{5.11}
\end{equation}

\begin{prop}\label{propmodstr5.2}
Let $U(t)$ satisfy \eqref{modstr.1} and \eqref{modstr.2}.  For any
$\gamma\ge 2\vee (2/\delta)$, we have
\begin{equation}
\|\mathfrak{U}f\|_{\ell^{q}_\Box (L^\infty(\mathbb{R}, L^2)) }
\le C  \| f\|_{\ell^{-\alpha/2, q}_\Box (L^{\gamma'}(\mathbb{R},
L^{p'}))}. \label{modstr.12}
\end{equation}
In addition, if $\gamma' \le q$, then
\begin{equation}
\|\mathfrak{U}f\|_{L^\infty(\mathbb{R}, M_{2,q})} \le C  \|
f\|_{L^{\gamma'}(\mathbb{R}, M^{-\alpha/2}_{p',q})}. \label{modstr.13}
\end{equation}
\end{prop}

\begin{prop}\label{propmodstr5.4}
Assume that $U(t)$ satisfies \eqref{modstr.1} and \eqref{modstr.2},
$\gamma\ge \max (2/\delta, 2)$. Then we have
\begin{equation}
\|\mathfrak{U}f\|_{\ell^{\alpha/2, q}_\Box (L^\gamma(\mathbb{R},
L^p)) } \le C  \| f\|_{\ell^q_\Box (L^1(\mathbb{R}, L^2))}.
\label{modstr.17}
\end{equation}
In addition, if  $\gamma\ge q$, then
\begin{equation}
 \|\mathfrak{U}f\|_{L^\gamma(\mathbb{R}, M^{\alpha/2}_{p,q})} \le C
\| f\|_{L^1(\mathbb{R}, M_{2, q})}. \label{modstr.18}
\end{equation}
\end{prop}

The Schr\"odinger semi-group corresponds to the cases $\alpha=0$,
$\delta=n(1/2-1/p)$ and $2\le p<\infty$. Taking $q=1$ in Propositions \ref{propmodstr5.1}--\ref{propmodstr5.4}, we immediately have

\begin{cor}\label{propmodstr5.5}
Let $2\le p< \infty$, $\gamma \ge 2\vee \gamma(p)$, and
\begin{equation}
 \frac{2}{\gamma(p)}=n\Big(\frac{1}{2}-\frac{1}{p}\Big). \label{modstr.20}
\end{equation}
Let $S(t)= e^{{\rm i}t\Delta}$, $\mathfrak{A}= \int^t_0 S(t-s)
\cdot ds$. Then
\begin{equation}
\left\|S(t) \varphi \right\|_{\ell^1_\Box (L^\gamma(\mathbb{R}, L^p)
)}   \le C  \| \varphi\|_{M_{2, 1}},  \label{modstr.21}
\end{equation}
\begin{equation}
\left\|\mathfrak{A} f \right\| _{\ell^1_\Box (L^\gamma (\mathbb{R},
L^p)) \cap \ell^1_\Box (L^\infty(\mathbb{R}, L^2)) }   \le C
\|f\|_{\ell^1_\Box (L^{\gamma'}(\mathbb{R}, L^{p'}))}. \label{modstr.22}
\end{equation}
\end{cor}

Similar to Corollary \ref{propmodstr5.5}, we have

\begin{cor}\label{propmodstr5.6}
Let $2\le p< \infty$, $\theta\in (0,1]$, $1\le q<\infty$,
\begin{equation}
 \frac{2}{\gamma_\theta (p)}=n\theta \Big(\frac{1}{2}-\frac{1}{p}\Big), \quad
 2\sigma=(n+2)\theta \Big(\frac{1}{2}-\frac{1}{p}\Big).  \label{modstr.24}
\end{equation}
Let $G(t)$ be as in \eqref{modecay.5}, $\mathfrak{G}=\int^t_0 G(t-s)
\cdot ds$. Then for any $\gamma \ge 2\vee \gamma_\theta (p)$, we have
\begin{equation}
\left\|G(t) \varphi \right\|_{\ell^{-\sigma, q}_\Box (L^\gamma(\mathbb{R}, L^p))}
 \le C  \| \varphi\|_{M_{2, q}},  \label{modstr.25}
 \end{equation}
 \begin{equation}
\left\|\mathfrak{G} f \right\| _{\ell^{-\sigma, q}_\Box
(L^\gamma(\mathbb{R}, L^p)) \cap \ell^q_\Box (L^\infty(\mathbb{R},
L^2))}   \le C  \|f\|_{\ell^{\sigma, q}_\Box
(L^{\gamma'}(\mathbb{R}, L^{p'}))}. \label{modstr.26}
\end{equation}
\end{cor}
Related Strichartz estimates in Wiener amalgam spaces for the Schr\"odinger equation were obtained by Cordero and  Nicola \cite{CorNik08}.

\subsection{Wellposedness for NLS and NLKG}

We study the Cauchy problem for NLS and give partial answers to Question 2.3. Let us consider
\begin{equation}
{\rm i}u_t + \Delta u = f(u), \quad u(0,x)= u_0(x). \label{NLSmod}
\end{equation}
Noticing that $B^n_{\infty,1} \subset M_{\infty, 1} \subset B^0_{\infty,1}\subset L^\infty$ are sharp embeddings, up to now, we can not get the wellposedness of NLS in  $L^\infty$ or in $B^0_{\infty,1}$.  However, we can obtain the local wellposedness of NLS in $M_{\infty, 1}$.
We have
\begin{thm} [\cite{BenOk08, CorNik09}]
\label{thmmodnls1.1}
Let $n\ge 1$, $f(u)= \lambda |u|^{\kappa}u$, $\kappa\in 2\mathbb{N}$, $\lambda \in \mathbb{R}$,
 $u_0 \in M_{p,1}$ and $1\le p \le \infty$. Then there exists a $T>0$ such that
\eqref{NLSmod} has a unique solution $u\in C([0,T), M_{p,1})$.  Moreover, if $T<\infty$, then $\lim\sup_{t\nearrow T}\|u(t)\|_{M_{p,1}}=\infty$.
\end{thm}

If the nonlinearity has an exponential growth, say $f(u)=\lambda(e^{|u|^2}-1)u$, the result in Theorem \ref{thmmodnls1.1} also holds.
Noticing that $B^{n/2}_{2,1} \subset M_{2, 1} \subset B^0_{2,1} \cap C(\mathbb{R}^n)$ are  sharp embeddings,
we can get that NLS is global wellposed in $M_{2,1}$ if initial data are sufficiently small.

\begin{thm} [\cite{WaHe07}] \label{thmmodnls1.3}
Let $n\ge 1$, $f(u)= \lambda |u|^{\kappa}u$, $\kappa\in 2\mathbb{N}$, $\lambda \in \mathbb{R}$,
$\kappa \ge 4/n$, $u_0 \in M_{2,1}$ and there exists a sufficiently small
  $\delta>0$ such that $\|u_0 \|_{ M_{2,1}} \le \delta$. Then
\eqref{NLSmod} has a unique solution
\begin{equation}
u\in C(\mathbb{R}, M_{2,1}) \cap \ell^1_\Box (L^{p}_{x, t\in
\mathbb{R}}), \label{modnls1.26}
\end{equation}
where $p\in [2+4/n, \, 2+\kappa] \cap \mathbb{N}$,
$\ell^1_\Box (L^{p}_{x, t\in \mathbb{R}}) $ is as in \eqref{fsud}.
\end{thm}

\begin{thm} [\cite{WaHe07}] \label{thmmodnls1.4}
Let $n\ge 2$, $f(u)= \lambda (e^{\varrho |u|^2}-1) u$,
$\lambda\in \mathbb{C}$ and $\varrho>0$. Assume that $u_0 \in M_{2,1}$
and there exists a sufficiently small $\delta>0$ such that $\|u_0 \|_{ M_{2,1}}
\le \delta$. Then \eqref{NLSmod} has a unique solution
\begin{equation}
u\in C(\mathbb{R}, M_{2,1}) \cap \ell^1_\Box (L^{4}_{x, t\in
\mathbb{R}}). \label{modnls1.28}
\end{equation}
\end{thm}

We now consider the initial value problem for NLKG,
\begin{equation}
  u_{tt} +(I-\Delta) u+f(u) =0, \quad u(0)=u_0,\; u_t(0)=u_1.
\label{nlkgmod1.2}
\end{equation}

Analogous to NLS, we have

\begin{thm} [\cite{BenOk08, CorNik09}] \label{thmmodnlkg1.2}
Let $n\ge 1$, $f(u)= \lambda |u|^{\kappa}u$, $\kappa\in 2\mathbb{N}$, $\lambda \in \mathbb{R}$,
 $(u_0, u_1) \in M_{p,1}\times M^{-1}_{p,1}$ and $1\le p \le \infty$. Then there exists a $T>0$ such that
\eqref{NLSmod} has a unique solution $(u, u_t)\in C([0,T), M_{p,1}) \times C([0,T), M^{-1}_{p,1})$.
Moreover, if  $T<\infty$, then $\lim\sup_{t\nearrow T}(\|u(t)\|_{M_{p,1}}+ \|u_t(t)\|_{M^{-1}_{p,1}})=\infty$.
\end{thm}

If the nonlinearity has an exponential growth, the corresponding results as in Theorem \ref{thmmodnlkg1.2} also hold.

\begin{thm} [\cite{WaHe07}] \label{thmmodnlkg1.6}
Let $n\ge 1$, $f(u)= \lambda u^{1+\kappa} $, $\kappa\in \mathbb{N}$ and
$\kappa \ge 4/n$. Put
\begin{equation}
\sigma= \frac{n+2}{n(2+\kappa)}. \label{modgwp1.31}
\end{equation}
Assume that $(u_0,u_1) \in M^\sigma_{2,1}\times M^{\sigma-1}_{2,1}$
and there exists a sufficiently small $\delta>0$ such that $\|u_0 \|_{
M^\sigma_{2,1}}+ \|u_1\|_{ M^{\sigma-1}_{2,1}} \le \delta$. Then
\eqref{nlkgmod1.2} has a unique solution
\begin{equation}
u\in C(\mathbb{R}, M^\sigma_{2,1}) \cap \ell^1_\Box
(L^{2+\kappa}_{x, t\in \mathbb{R}}). \label{modgwp1.32}
\end{equation}
\end{thm}

\begin{thm} [\cite{WaHe07}] \label{thmmodnlkg1.7}
Let $n\ge 2$, $f(u)= \sinh u -u$ and $\sigma=(n+2)/4n$. Assume that
$(u_0,u_1) \in M^\sigma_{2,1}\times M^{\sigma-1}_{2,1}$ and there exists a sufficiently small $\delta>0$ such that $\|u_0 \|_{ M^\sigma_{2,1}}+
\|u_1\|_{ M^{\sigma-1}_{2,1}} \le \delta$. Then \eqref{nlkgmod1.2} has a unique solution
\begin{equation}
u\in C(\mathbb{R}, M^\sigma_{2,1}) \cap \ell^1_\Box (L^{4}_{x, t\in
\mathbb{R}}). \label{modgwp1.34}
\end{equation}
\end{thm}

\section{Derivative nonlinear Schr\"odinger equations}

We study the initial value problem for the derivative nonlinear Schr\"odinger equation (gDNLS)
\begin{equation}
{\rm i} u_t + \Delta_\pm u = F(u, \bar{u}, \nabla u, \nabla
\bar{u}), \quad u(0,x)= u_0(x),
 \label{gNLS}
\end{equation}
where $u$ is a complex-valued function of $(t,x)\in \mathbb{R}
\times \mathbb{R}^n$,
\begin{equation}
\Delta_\pm u = \sum^n_{i=1} \varepsilon_i \partial^2_{x_i}, \quad
\varepsilon_i \in \{1,\,  -1\}, \quad i=1,...,n,
 \label{Delta-pm}
\end{equation}
$\nabla =(\partial_{x_1},..., \partial_{x_n})$, $F:
\mathbb{C}^{2n+2} \to \mathbb{C}$ is a series of $z\in \mathbb{C}^{2n+2} $,
\begin{equation}
F(z) = F(z_1,..., z_{2n+2})= \sum_{3 \le |\beta|<\infty} c_\beta
z^\beta, \quad c_\beta \in \mathbb{C},
 \label{poly}
\end{equation}
$|c_\beta |\le C^{|\beta|}$. The typical nonlinear term  is
$$
F(u, \bar{u}, \nabla u, \nabla \bar{u}) = |u|^2\vec{\lambda}\cdot
\nabla u + u^2\vec{\mu}\cdot \nabla \bar{u} +|u|^2 u,
$$
see \cite{ClTu90,DT89,TD89}. Another model is
$$
F(u, \bar{u}, \nabla u, \nabla \bar{u})= (1+ |u|^2)^{-1}(\nabla
u)^2 \bar{u} = \sum^\infty_{k=0} (-1)^k |u|^{2k} (\nabla u)^2 \bar{u},  \quad |u|<1,
$$
which is an equivalent version of the Schr\"odinger flow \cite{Io-Ke07}. The non-elliptic gDNLS arises  in the strongly interacting many-body
systems near the criticality, where anisotropic interactions are manifested by the presence of
the non-elliptic case, as well as additional residual terms which
involve cross derivatives of the independent variables \cite{ClTu90, DT89, TD89}. Some water wave  and completely integrable system models in higher spatial dimensions are also non-elliptic, cf. \cite{AbHa75, ZaKu86, ZaSc80}.  A large amount of work has been devoted to the study of gDNLS, see \cite{HaOzSJMA94, Io-Ke07, KePoRoVe, KePoVe1, KePoVe2, Klai, Kl-Po, Oz-Zh08, RuSu2, Sh82, SuRu05}.

Since the nonlinearity in gDNLS contains derivative terms and the Strichartz inequalities can not absorb any derivatives, gDNLS can not be solved  if we use only the Strichartz estimate. One needs to look for some other ways to handle the derivative terms in the nonlinearity. Up to now, three kinds of methods seem to be very useful for gDNLS. One is to use the energy estimate to deal with the derivatives in the nonlinearity, the second way is to use Bourgain's  space $X^{s,b}$ and the third technique is Kato's smooth effect estimates.
Of course, there are some connections between these methods.

In this survey paper we only discuss the smooth effect estimates together with the frequency-uniform decomposition techniques for gDNLS and we show that it is globally wellposed and scattering in a class of modulation spaces.

For convenience, we denote
$$
S(t)= e^{{\rm i}t \Delta_\pm}= \mathfrak{F}^{-1} e^{{\rm i}t
\sum^n_{j=1}\varepsilon_j \xi^2_j}\mathfrak{F}, \ \ \mathfrak{A} f
(t,x)= \int^t_0 S(t-\tau) f(\tau,x) d\tau.
$$

\medskip

We now state a global wellposedness and scattering result for gDNLS in modulation spaces.  We
denote by $L^{p_1}_{x_i} L^{p_2}_{(x_j)_{j\not=i}} L^{p_2}_t:=
L^{p_1}_{x_i} L^{p_2}_{(x_j)_{j\not=i}} L^{p_2}_t(\mathbb{R}^{1+n})$
the anisotropic Lebesgue space for which the norm is defined by
\begin{align}
\|f\|_{L^{p_1}_{x_i} L^{p_2}_{(x_j)_{j\not=i}}L^{p_2}_t} =\left\|
\|f\|_{L^{p_2}_{x_1,...,x_{j-1}, x_{j+1},...,x_n}
L^{p_2}_t(\mathbb{R} \times \mathbb{R}^{n-1})}
\right\|_{L^{p_1}_{x_i}(\mathbb{R})} . \label{Notation.1}
\end{align}
For $k=(k_1,...,k_n)$,  we write
\begin{eqnarray}
\|u\|_{X^s_\alpha} =    \ \sum_{i, \, \ell=1}^n \ \sum_{k\in
\mathbb{Z}^n, \ |k_i|>10 \vee \max_{j\neq i} |k_j|} \langle k_i\rangle^{s-1/2}
\left\|\partial^\alpha_{x_\ell} \Box_k u \right\|_{L^\infty_{x_i}
L^2_{(x_j)_{j\not= i}}L^2_t } \nonumber\\
   +  \ \sum_{i,\, \ell=1}^n \ \sum_{k\in \mathbb{Z}^n}  \left\|
\partial^\alpha_{x_\ell} \Box_k u \right\|_{L^{2}_{x_i}
L^\infty_{(x_j)_{j\not= i} }
L^\infty_t },   \\
\|u\|_{S^{s}_\alpha}   =    \ \sum_{\ell=1}^n \sum_{k\in
\mathbb{Z}^n} \langle k \rangle^{s-1} \left\|
\partial^\alpha_{x_\ell} \Box_k u \right\|_{L^\infty_t L^2_x  \bigcap  L^3_t L^{6}_{x}
 },   \\
 \|u\|_{X^s}=
\sum_{\alpha=0,1}  \|u\|_{X^s_\alpha},   \quad
 \|u\|_{S^s}=  \sum_{\alpha=0,1} \|u\|_{S^s_\alpha}.
 \label{sum}
 \end{eqnarray}

\begin{thm} [\cite{WaHaHu09, Wa11}] \label{DNLS-modm2}
Let $n\ge 3$,  $u_0 \in M^{3/2}_{2,1}$ and there exists a suitably small $\delta>0$ such that
$\|u_0\|_{M^{3/2}_{2,1}} \le \delta$. Then
\eqref{gNLS} has a unique solution $u\in C(\mathbb{R}, M^{3/2}_{2,1}) \cap X^{3/2} \cap S^{3/2}$,
$\|u\|_{X^{3/2} \cap S^{3/2}} \le C  \delta$. Moreover, the scattering operator $S$ of
 \eqref{gNLS} carries a whole zero neighborhood in $C(\mathbb{R}, M^{3/2}_{2,1})$ into $C(\mathbb{R},
M^{3/2}_{2,1})$.
\end{thm}

\begin{rem}
Recently, this technique was also developed for the Navier-Stokes equation and the dissipative nonlinear electrohydrodynamic system \cite{Wang8, Iw10, DeZhCu10}.
\end{rem}

\section{Canonical transformations}

It is interesting to generalize results in previous two sections
to the case of dispersive operators $a(D)$ instead of $-\Delta$.
Recently Ruzhansky and Sugimoto has introduced a new
idea to establish fundamental estimates for dispersive equations
based on the idea of canonical transformations, and this attempt
is quite successful for smoothing estimates (\cite{RuSu1}).

Let $\psi:\mathbb{R}^n\setminus0\to\mathbb{R}^n\setminus0$ be $C^\infty$-maps
satisfying $\psi(\lambda\xi)=\lambda\psi(\xi)$
for all $\lambda>0$ and $\xi\in\mathbb{R}^n\setminus0$, and let
\begin{equation}\label{DefofI}
\begin{aligned}
 I_\psi u(x)&=\mathfrak{F}^{-1}[(\mathfrak{F}u)(\psi(\xi))](x)
\\
&=(2\pi)^{-n}\int_{\mathbb{R}^n}\int_{\mathbb{R}^n}
      e^{{\rm i}(x\cdot\xi-y\cdot\psi(\xi))}u(y)\, dy d\xi.
\end{aligned}
\end{equation}
We remark that we have the formula
\begin{equation}\label{eq:cnon}
a(D)\cdot I_\psi=I_\psi\cdot\sigma(D),\quad
a(\xi)=(\sigma\circ\psi)(\xi).
\end{equation}
For example, for a positive function $a(\xi)$
satisfying $a(\lambda\xi)=\lambda^2 a(\xi)$
for all $\lambda>0$ and $\xi\in\mathbb{R}^n\setminus0$,
we have
\[
a(D)\cdot I_\psi=I_\psi\cdot(-\Delta)
\]
if we take
\[
\sigma(\eta)=|\eta|^2,\quad
\psi(\xi)=\sqrt{a(\xi)}\frac{\nabla a(\xi)}{|\nabla a(\xi)|},
\]
provided that $\nabla a(\xi)\not=0$ for $\xi\not=0$.
The latter is achieved if we assume that the Gaussian curvature of the
hypersurface $\Sigma=\{\xi:\ a(\xi)=1\}$ never vanishes,
in which case also $I_\psi$ has the inverse $I_{\psi^{-1}}=I_\psi^{-1}$ because the Gauss map
$\Sigma\ni\xi\mapsto\frac{\nabla a(\xi)}{|\nabla a(\xi)|}\in {\mathbb S}^{n-1}$
is a diffeomorphism.

Thus we can induce the same estimate
for $a(D)$ from the estimates for $-\Delta$
if we establish the boundedness of operators $I_\psi$ and $I_{\psi^{-1}}$
on modulation spaces.
Ruzhansky, Sugimoto, Toft and Tomita \cite{RuSuTofTom} discuss
such boundedness properties, and we have a positive result for the local boundedness.
Let $L^q_s$ be the space of functions such that
$\langle x\rangle^s f\in L^q$:

\begin{thm}\label{TH:local-lq}
Let $s\in\mathbb{R}$, $1 \le p,q \le \infty$,
and let $\psi : \mathbb{R}^n \to \mathbb{R}^n$ be such that
the pullback $\psi^*:f\mapsto f\circ \psi$ is bounded
on $L^q_s(\mathbb{R}^n)$.
Then
$I_\psi$ is locally bounded on $M_{p,q}^s(\mathbb{R}^n)$.
\end{thm}

By Theorem \ref{TH:local-lq} and Propositions
\ref{propmodpp.1}, \ref{propmodecay.1}, we straightforwardly obtain estimates
\eqref{modpp.5} and \eqref{modecay.4}
for $S(t)=\chi(x)e^{-{\rm i}ta(D)}\tilde\chi(x)$
(and we have to also restrict the case to $1\le q \le \infty$)
form the estimates for $S(t)=e^{{\rm i}t\Delta}$,
where $\chi,\tilde\chi$ are cut-off functions.

As for the global boundedness of $I_{\varphi}$, we have unfortunately
a negative result:
\begin{thm}\label{TH:B-H}
Let $1 \le p,q \le \infty$, $ 2 \neq p<\infty$,
and let $\psi : \mathbb{R}^n \to \mathbb{R}^n$ be a $C^1$-function.
Assume that operator $I_\psi$
is bounded on $M_{p,q}(\mathbb{R}^n)$.
Then $\psi$ is an affine mapping.
\end{thm}

\section{Open questions}

It is known that the algebra property of function spaces is of importance for PDE. Up to now the following question is not clear for us:

\begin{que}
Let $\alpha \in (0, \infty)$. Does
$$
\||u|^\alpha u\|_{M_{p,1}} \le C \| u\|^{\alpha+1}_{M_{p,1}}
$$
hold for all $u\in M_{p,1}$?
\end{que}
It is known that if $\alpha\in 2\mathbb{N}$, the answer is affirmative. If $\alpha$ is not an even integer, the question seems very difficult, cf. \cite{SuToWa10}.

The global well posedness of NLS in modulation spaces for large initial data seems open.

\begin{que}
Can we show that NLS \eqref{NLSmoda} is global well posed  if the initial data $u_0\in M_{2,1}$ is large?
\end{que}

We should also discuss non-affine transforms which induce the globally bounded
canonical transformations.
Note that such transforms must not be $C^1$-mappings in view of
Theorem \ref{TH:B-H}.

\begin{que}
Can we show the global boundedness of the operator $I_\psi$ in \eqref{DefofI}
on modulation spaces for a homogeneous change of variables $\psi$?
\end{que}
Some partial answers to this question appeared in \cite{RuSuTofTom}.

% ------------------------------------------------------------------------

% ------------------------------------------------------------------------
\end{document}